\documentclass[amssymb,11pt]{amsart}
\usepackage{latexsym}

\newlength{\standardunitlength}
\newtheorem{prop}{Proposition}[section]

\newtheorem{lemma}[prop]{Lemma}
\newtheorem{cor}[prop]{Corollary}
\newtheorem{theorem}[prop]{Theorem}

\begin{document}

\begin{center} {\bf A Card Shuffling Analysis of Deformations of the Plancherel
Measure of the Symmetric Group}
\end{center}

\begin{center}
By Jason Fulman
\end{center}

\begin{center}
University of Pittsburgh
\end{center}

\begin{center}
Department of Mathematics
\end{center}

\begin{center}
301 Thackeray Hall
\end{center}

\begin{center}
Pittsburgh, PA 15260
\end{center}

\begin{center}
Email: fulman@math.pitt.edu
\end{center}

\begin{center}
Version of 2/20/03
\end{center}

\newpage

{\bf Abstract}: We study deformations of the Plancherel measure of the symmetric group by lifting them to the symmetric group and using combinatorics of card shuffling. The existing methods for analyzing deformations of Plancherel measure are not obviously applicable to the examples in this paper. The main idea of this paper is to find and analyze a formula for the total variation distance between iterations of riffle shuffles and iterations of ``cut and then riffle shuffle''. Similar results are given for affine shuffles, which allow us to determine their convergence rate to randomness. 

\section{Introduction}

	In recent years there has been significant interest in the
	study of longest increasing subsequences in random
	permutations. It is beyond the scope of this paper to survey
	the subject, but the connections are fascinating and relate to
	random matrix theory, Painleve functions, Riemann surfaces,
	point processes, Riemann Hilbert problems, and much
	else. Recent surveys include \cite{AD} and \cite{De}.

	By means of the Robinson-Schensted-Knuth (RSK) correspondence, the uniform measure on the symmetric group induces a measure (called the Plancherel measure) on integer partitions of the number $n$. We assume that the reader is familiar with these concepts; a very readable treatment is the survey \cite{AD}. One of the more interesting aspects of the theory is that to analyze the Plancherel measure, one often analyzes a deformation of the Plancherel measure and then takes a limit in which the deformation approaches Plancherel measure. This approach is used, for instance, in \cite{J}; see also the survey \cite{BO}.

	Next let us recall the k-riffle shuffle measure on a deck of $n$ cards, defined in \cite{BD} as a generalization of the Gilbert-Shannon-Reeds model of the way real people shuffle cards. For an integer $k \geq 1$, this shuffle may be described as follows. Given \
a deck of $n$ cards, one cuts it into $k$ piles with probability of pile sizes \
$j_1,\cdots,j_k$ given by $\frac{{n \choose j_1,\cdots,j_k}}{k^n}$. Then cards \
are dropped from the packets with probability proportional to the pile size at \
a given time. Thus if the current pile sizes are $A_1,\cdots,A_k$, the next card is dropped from pile $i$ with probability $A_i/(A_1+\cdots+A_k)$. Let $d(\pi)$ be the number of descents of $\pi$ (i.e. the number of $i$ with $1 \leq i \leq n-1$ such that $\pi(i) > \pi(i+1)$). There is a simple formula for the chance that a $k$-shuffle gives the inverse of a permutation $\pi$ in terms of $d(\pi)$. Denoting this probability by $R_{k,n}(\pi)$, it is proved in \cite{BD} \
that \[ R_{k,n}(\pi) = \frac{{n+k-d(\pi)-1 \choose n}}{k^n}.\] 

	One immediate observation is that as $k \rightarrow \infty$, the measure $R_{k,n}$ converges to the uniform distribution. In fact, much more is true. The paper \cite{BD} establishes the following equation in the group algebra of the symmetric group, for any integers $k_1,k_2 \geq 1$: \[ \left( \sum_{\pi \in S_n} R_{k_1,n}(\pi) \pi \right) \left( \sum_{\pi \in S_n} R_{k_2,n}(\pi) \pi \right) = \left( \sum_{\pi \in S_n} R_{k_1k_2,n}(\pi) \pi \right).\] Thus the riffle shuffle measure $R_{k,n}$ is {\it closed under iteration}, in the sense that a $k_1$ shuffle followed by a $k_2$ shuffle is equal to a $k_1k_2$ shuffle. So iterating a $2$-shuffle $k$ times is equivalent to performing a single $2^k$ shuffle, and one can study properties of iterates of riffle shuffles using the explicit formula for $R_{k,n}$.

	To make this quantitative, one often uses the total variation distance $||P_1-P_2||$ between two probability distributions $P_1$ and $P_2$ on a finite set $X$. It is defined as \[ \frac{1}{2} \sum_{x \in X} |P_1(x)-P_2(x)|.\] One reason that this is a useful method of quantifying the distance between two distributions is that the total variation distance between probability distributions $P_1,P_2$ can also be defined as \[ max_{S \subset X} |P_1(S)-P_2(S)|.\] Thus when the total variation distance is small, any event has nearly the same probability under the measures $P_1,P_2$. Aldous \cite{A} showed that for large $n$, $\frac{3}{2} log_2(n)$ 2-riffle shuffles are necessary and suffice to randomize a deck of $n$ cards. This was sharpened by Bayer and Diaconis \cite{BD} who showed that (for any $n$), if $k = \frac{3}{2} \log_2 n + c$ and $U_n$ denotes the uniform distribution on $S_n$, then $|| R_{2^k,n} - U_n ||$ is
\begin{equation*}
1 - 2 \Phi \Bigg(\frac{-2^{-c}}{4 \sqrt3}
\Bigg) + 0 \Bigg( \frac{1}{n^{1/4}} \Bigg) \quad \mbox{where} \quad
\Phi(x) = \frac{1}{\sqrt{2 \pi}} \int_{\infty}^x e^{-t^2 / 2} dt .
\end{equation*} This marvelous result (sometimes loosely stated as $\frac{3}{2}log_2(n)$ iterations of a 2-riffle shuffle are necessary and suffice for randomness) is also quite powerful. As a consequence of the Bayer/Diaconis result, it follows that after $\frac{3}{2}log_2(n)$ iterations of a 2-riffle shuffle, the distribution of the longest increasing subsequence is close to that of a random permutation.

	The reader might at this point wonder what card shuffling has
	to do with deformations of the Plancherel measure of the
	symmetric group. The point is a result of Stanley \cite{S}
	that the the distribution of the longest increasing
	subsequence (or more generally the RSK shape) under the
	measure $R_{k,n}$ is equal to the distribution of the longest
	weakly increasing subsequence in a random length $n$ word on
	$k$ letters. In fact Johansson's proof (\cite{J}) of the
	Baik-Deift-Johansson theorem used discrete orthogonal
	polynomial ensembles to analyze the longest weakly increasing
	subsequence in random length $n$ words on $k$ letters, and
	then let $k \rightarrow \infty$. As observed in \cite{Finc},
	Johansson's work implies that iterating a 2-riffle shuffle
	$\frac{5}{6}log_2(n)$ times brings the distribution of the
	longest increasing subsequence close to that of a random
	permutation. Note that this is sharper than the
	$\frac{3}{2}log_2(n)$ bound of the previous paragraph. Section
	\ref{shuffleinc} of this paper briefly summarizes these
	developments, giving precise statements and pointers to the
	literature.
		
	Given the above discussion, it is very natural to study
	deformations of Plancherel measure given by applying the RSK
	correspondence to other shuffles which are closed under
	iteration. The first question is: which shuffles are the most
	natural? The series of papers \cite{F1},\cite{F2},\cite{F3}
	used Lie theory to define a number of methods of shuffling
	which are closed under iteration and mathematically very
	natural. The paper \cite{Finc} showed that many, but not all,
	of these shuffles interact nicely with RSK correspondence and with Toeplitz determinants. This paper aims to analyze
	deformations of Plancherel measure arising from shuffles which on one hand are very natural, but whose interaction with the RSK correspondence is at the present writing unclear (essentially because the interaction of the cyclic descent set of a permutation with the RSK correspondence is unclear). 

	The first shuffle we analyze can be physically described as
	cutting the deck at a uniformly chosen random position and
	then performing a k-riffle shuffle. This shuffle was first
	studied in the paper \cite{F2}, although throughout that paper
	we incorrectly called it a shuffle followed by a cut instead
	of a cut followed by a shuffle; shuffles followed by cuts had
	been studied in \cite{BD}. In any case, the measure defined in
	\cite{F2} is that which picks a permutation $\pi$ on n-symbols
	with probability \[ C_{k,n}(\pi) = \frac{{n+k-cd(\pi)-1 \choose n-1}}{nk^{n-1}},\] where $cd(\pi)$ is defined as
	follows. A permutation $\pi$ on $n$ symbols is said to have a
	cyclic descent at position $n$ if $\pi(n)>\pi(1)$. Then
	$cd(\pi)$ is defined as $d(\pi)$ if $\pi$ has no cyclic
	descent at $n$ and as $d(\pi)+1$ if $\pi$ has a cyclic descent
	at $n$. Thus $C_{k,n}(\pi)$ is the chance of obtaining
	$\pi^{-1}$ after applying a cut and then a k-riffle
	shuffle. This probability measure
	shares many of the truly remarkable properties of the measure
	$R_{k,n}$. For example, the measure $C_{k,n}$, like the
	measure $R_{k,n}$, has the property that one can
	exactly understand its distribution on conjugacy classes, even
	though $C_{k,n}$ is not constant on conjugacy classes. Another
	property which is more relevant to the current discussion is
	that \[ \left( \sum_{w \in S_n} C_{k_1,n}(w) w \right) \left(
	\sum_{w \in S_n} C_{k_2,n}(w) w \right) = \left( \sum_{w \in
	S_n} C_{k_1k_2,n}(w) w \right).\] In words, if one iterates
	the process ``cut and then perform a k-riffle shuffle'', one
	obtains exactly the same result as from first performing a cut
	and then iterating the k-riffle shuffles. Thus only the
	initial cut contributes to the randomness. It follows that
	$\frac{3}{2}log_2(n)$ iterations of cut and then $2$-riffle
	shuffle are necessary and suffice to get close to
	randomness (see \cite{F2} for details). In particular, the longest
	increasing subsequence must be random after
	$\frac{3}{2}log_2(n)$ iterations.

	In fact since $C_{k,n}$ is the convolution of $R_{k,n}$ with the probability measure corresponding to cutting the deck at a random position, we have for free that $||C_{k,n}-U_n|| \leq ||R_{k,n}-U_n||$. Indeed, if $P,Q$ are any measures whatsoever on the
symmetric group, then the chance that $P$ followed by $Q$ yields a permutation $\pi$ is $\sum_{\tau} Q(\pi \tau^{-1}) P(\tau)$, and it is easy to see that \[ \frac{1}{2} \sum_{\pi} | \sum_{\tau} Q(\pi
\tau^{-1}) P(\tau) - \frac{1}{n!}| \leq ||Q-U_n|| \] and
\[ \frac{1}{2} \sum_{\pi} | \sum_{\tau} Q(\pi
\tau^{-1}) P(\tau) - \frac{1}{n!}| \leq ||P-U_n||. \] But although the
measure $C_{k,n}$ is closer to the uniform distribution than the
measure $R_{k,n}$, it does not follow that any statistic is more
random under the measure $C_{k,n}$ than under the measure
$R_{k,n}$. For example let $\pi$ be a permutation on 3 symbols which
has two cyclic descents but only 1 descent (for instance the
permutation transposing symbols 1 and 2). Then $C_{k,3}$ assigns
probability $\frac{1}{6}-\frac{1}{6k}$ to $\pi$, but the measure
$R_{k,3}$ assigns probability $\frac{1}{6}-\frac{1}{6k^2}$ to
$\pi$. Thus if $S$ is the event ``being equal to the transposition
(12)'', $|C_{k,3}(S)-U_3(S)| \geq |R_{k,3}(S)-U_3(S)|$. So although it
is probably true, we do not know that $\frac{5}{6}log_2(n)$ iterations
of a cut followed by a $2$ riffle shuffle suffice to randomize the
longest increasing subsequence.

	A main result of this paper is to derive sharp bounds on the total variation distance $||C_{k,n}-R_{k,n}||$. A wonderful simplification occurs and the absolute values signs in the definition of total variation distance evaporate and exact generating functions can be used (this does not happen, for example, in the Bayer/Diaconis analysis of $||R_{k,n}-U_n||$, where some serious effort is required). We conclude that whereas both a 2-riffle shuffle and a cut followed by a 2-riffle shuffle take $\frac{3}{2}log_2(n)$ steps to converge to the uniform distribution, only $log_2(n)$ steps are necessary for them to converge to each other. As a corollary, since \[ ||C_{k,n}-U_n|| \leq ||C_{k,n}-R_{k,n}|| + ||R_{k,n}-U_n|| \] we obtain that for a cut followed by a 2-riffle shuffle, $log_2(n)$ steps suffice to bring the distribution of the longest increasing subsequence close to that of a random permutation. This is a significant improvement to the $\frac{3}{2}log_2(n)$ bound from two paragraphs ago. And even if it turns out that one can get sharper bounds for the longest increasing subsequence under the measure $C_{k,n}$ by other methods, our method is applicable to other statistics as well. For instance since a cut followed by a riffle shuffle makes the position of any particular card perfectly random, we recover a recent result of \cite{U} stating that $log_2(n)$ steps suffice for a $2$-riffle shuffle to randomize any particular card. 

	The final section of this paper treats another variation of riffle shuffles called type $A$ affine k-shuffles \cite{F2},\cite{F3}. These are also closed under iteration, and the chance $A_{k,n}(\pi)$ of a permutation depends both on the number of cyclic descents and the major index of the permutation. It is proved that $log_2(n)$ steps suffice to bring a type $A$ affine 2-shuffle and a cut followed by a 2-riffle shuffle close to each other. This implies that $\frac{3}{2}log_2(n)$ steps are necessary and suffice for convergence of a type $A$ affine 2-shuffles to the uniform distribution, and that the longest increasing subsequence converges is close to that of a random permutation in at most $log_2(n)$ steps. 

	The precise organization of this paper is as follows. Section
	\ref{shuffleinc} very briefly reviews connections between
	riffle shuffles, longest increasing subsequences, and random
	matrices; some of these will be needed for this paper. Section
	\ref{totvarcut} derives an
	exact expression for $||C_{k,n}-R_{k,n}||$. Then it analyzes
	the formula and records some corollaries. Section \ref{affine}
	treats the measure $A_{k,n}$, proving similar results.

\section{Shuffling and Subsequences} \label{shuffleinc}

	This section very briefly recalls a few connections between shuffling, increasing subsequences, and random matrices. All of this material is known but perhaps not collected in one place, and we shall need some of it. 

	Let $L(\pi)$ denote the length of the longest increasing subsequence of a permutation chosen uniformly at random from the symmetric group on $n$ symbols. The Tracy-Widom distribution $F(t)$ is defined as \[ F(t)=exp \left(- \int_{t}^{\infty} (x-t)u^2(x)dx \right) \] where $u(x)$ is the solution of the Painleve II equation $u_{xx}=2u^3+xu$ such that $u$ is asymptotic to the negative of the Airy function as $x \rightarrow \infty$. The Baik-Deift-Johansson theorem (\cite{BDJ}), states that \[ lim_{n \rightarrow \infty} Prob \left(\frac{L-2 \sqrt{N}}{N^{1/6}} \leq t \right) = F(t) \] for all $t \in R$.
	
	Johansson \cite{J} later proved this using discrete orthogonal polynomial ensembles. He viewed the Poissonized Plancherel measure (the term Poissonized refers to the fact that the number of symbols $n$ is a Poisson random variable) as a limiting case of the Charlier ensemble, and showed that the Charlier ensemble could be exactly analyzed. The Charlier ensemble converges to the Poissonized Plancherel measure because random length $n$ words on $k$ symbols converge to random permutations on $n$ symbols as $k \rightarrow \infty$, and the Charlier ensemble is essentially what arises when applying the RSK to random words. Stanley \cite{S} observed that the distribution of the longest weakly increasing subsequence in a random length $n$ word on $k$ symbols has exactly the same distribution as the longest increasing subsequence in a random permutation after a $k$-riffle shuffle. This is simple to see. For example when $k=2$, given a word on the symbols $1,2$ such as \[ 1 \ 2 \ 1 \ 1 \ 2 \ 1 \ 1 \ 2 \ 1 \] one obtains a shuffle of the stacks $1,2,3,4,5,6$ and $7,8,9$ by putting the numbers in the first stack in the positions of the ones (going from left to right), then putting then numbers in the second stack in the position of the twos, and so on. In this example the permutation which emerges is \[ 1 \ 7 \ 2 \ 3 \ 8 \ 4 \ 5 \ 9 \ 6.\] Clearly a weakly increasing subsequence in the word corresponds to an increasing subsequence in the permutation associated to the word. Thus the study of longest increasing subsequences in random words (and hence the Charlier ensemble) is really the study of the RSK shape obtained after a permutation distributed as a riffle shuffle.

	In fact Johansson \cite{J} used orthogonal polynomial ensembles to prove that if $L(w)$ denotes the length of the longest weakly increasing subsequence in a random length $n$ word on $k$ symbols, and if $n,k \rightarrow \infty$ in such a way that $(log(n))^{1/6}/k \rightarrow 0$, then \[ lim_{n \rightarrow \infty} Prob \left( \frac{L(w)-\frac{n}{k}-2 \sqrt{n}}{\left( 1+\frac{\sqrt{n}}{k} \right)^{2/3} n^{1/6}} \leq t \right)= F(t) \] where as above $F(t)$ is the Tracy-Widom distribution. The paper \cite{Finc} translated this into card-shuffling language by noting that if $k = \lceil e^cn^{5/6} \rceil$, then \[ lim_{n \rightarrow \infty} Prob \left( \frac{L(w)-2 \sqrt{n}}{n^{1/6}} \leq t \right) = F(t-e^{-c}). \] Thus $\frac{5}{6}log_2(n)$ 2-riffle shuffles are necessary and suffice to bring the longest increasing subsequence close to that of a random permutation.

	We remark that there are connections with random matrices. First, it is known that the Tracy-Widom distribution is the limit distribution for the largest eigenvalue of a random GUE (Gaussian Unitary Ensemble) matrix \cite{TW1}. Moreover, the random word problem (or equivalently riffle shuffling) is related both to the ensemble of traceless GUE matrices and to the Laguerre ensemble (see Theorems 3 and 4 of \cite{TW2}, \cite{K} and the references in these papers for details).

\section{Shuffles with Cuts} \label{totvarcut}

	Throughout this section, we adhere to the notation in the introduction. Thus $d(\pi)$ and $cd(\pi)$ denote the number of descents and cyclic descents of a permutation. We also freely use the explicit formulas for $R_{k,n}$ and $C_{k,n}$ mentioned in the introduction.

	Although we will not use it in this paper, we record the following observation.

\begin{prop} A k-riffle shuffle followed by a cut and a cut followed by a k-riffle shuffle induce the same measure on conjugacy classes. \end{prop}

\begin{proof} Let $\zeta$ denote the cycle $(1 \mapsto 2 \cdots \mapsto n \mapsto 1)$. Let $s_k$ be defined as the element of the group algebra of $S_n$ equal to $\sum_{\pi \in S_n} R_{k,n}(\pi^{-1}) \pi$. Then a k-riffle shuffle followed by a cut is given by $s_k \left(\frac{1}{n}
 \sum_{j=1}^n \zeta^j \right)$ and a cut followed by a k-riffle
 shuffle is given by $\left( \frac{1}{n} \sum_{j=1}^n \zeta^j \right) s_k$. But it
 is straightforward to see that if $g$ is any element of the group
 algebra of the symmetric group, then $g \left( \sum_{j=1}^n \zeta^j
 \right)$ and $\left( \sum_{j=1}^n \zeta^j \right) g$ induce the same
 measure on conjugacy classes. Indeed, $g \zeta^j$ is conjugate to $\zeta^{n-j} g$. \end{proof}

	In what follows, we let $A_{n,d}$ denote the number of
	permutations on $n$ symbols with $d-1$ descents. Lemma
	\ref{Worpitzky} collects known properties of Eulerian numbers;
	see for example \cite{FoS}. Part 3 is called Worpitzky's
	identity and is equivalent to the fact that $R_{k,n}$ is a
	probability measure.

\begin{lemma} \label{Worpitzky} Let $A_{n,d}$ denote the number of
	permutations on $n$ symbols with $d-1$ descents.
\begin{enumerate}
\item \[ A_{n,d}=A_{n,n+1-d}.\]
\item \[ \sum_{\pi \in S_n} \frac{t^{d(\pi)+1}}{(1-t)^{n+1}} = \sum_{k \geq 0} k^n t^k.\]
\item \[ k^n = \sum_{1 \leq j \leq n} A_{n,j} {k+n-j \choose n}.\]
\end{enumerate}
\end{lemma}

\begin{lemma} \label{cyclicdes} For $n>1$, the number of permutations on $n$ symbols with $d$ cyclic descents and $d-1$ descents is $dA_{n-1,d}$, and the number of permutations on $n$ symbols with $d$ cyclic descents and $d$ descents is $(n-d)A_{n-1},d$. 
\end{lemma}

\begin{proof} It is proved in \cite{F2} that for $n>1$, the number of permutations on $n$ symbols with $d$ cyclic descents is equal to $nA_{n-1,d}$. Now any the set of permutations on n symbols with $d$ cyclic descents falls into equivalence classes under the operation of cyclic shifts. For example the permutation $1 3 5 2 4$ is a member of a set of 5 permutations \[ 1 3 5 2 4 , 3 5 2 4 1 , 5 2 4 1 3 , 2 4 1 3 5 , 4 1 3 5 2 .\] The number of cyclic descents is constant among this class, but clearly the proportion of members of this class with $d-1$ descents is $\frac{d}{n}$ and the proportion of members of this class with $d$ descents is $\frac{n-d}{n}$. \end{proof}

	Theorem \ref{exactformula} gives three expressions for the total variation distance between the measures $R_{k,n}$ and $C_{k,n}$. It is rather remarkable that such a simple formula (i.e. one not involving absolute values) exists. We remark that for the purpose of an upper bound, only the first equality in Theorem \ref{exactformula} is used. The third equality is useful for lower bounds and involves $B_n$, the nth Bernoulli number. They are defined by the generating function $\sum_{n \geq 0} \frac{B_nz^n}{n!} = \frac{z}{e^z-1}$, so that $B_0=1,B_1=-\frac{1}{2},B_2=\frac{1}{6},B_3=0,B_4=-\frac{1}{30}$ and $B_i=0$ if $i \geq 3$ is odd.  Also $(n)_t$ denotes $n(n-1)\cdots(n-t+1)$.

\begin{theorem} \label{exactformula} \begin{eqnarray*}
||R_{k,n}-C_{k,n}|| & = & \frac{1}{k^n n} \sum_{j=1}^{n-1} j(n-j) A_{n-1,j} {n+k-j-1 \choose n-1}\\
& = & \frac{1}{k^n} \left(k^n-\frac{k^{n+1}}{n}+(n+1) \sum_{s=1}^{k-1}s^n - k(n-1) \sum_{s=1}^{k-1} s^{n-1} \right)\\
& = & \left( \sum_{t=1}^{n-2} \frac{B_{t+1}}{k^t} \left( \frac{(n)_t}{t!} + \frac{(n-1)_t}{(t+1)!} \right) \right) + \frac{(n+1) B_n}{k^{n-1}}. \end{eqnarray*} \end{theorem}

\begin{proof} To compute $||R_{k,n}-C_{k,n}||$, note that for any
permutation $\pi$, either $cd(\pi)=d(\pi)+1$ or $cd(\pi)=d(\pi)$. If
$cd(\pi)=d(\pi)+1$, then by the explicit formulas for $R_{k,n}$ and
$C_{k,n}$ in the introduction, it follows that
\begin{eqnarray*} |R_{k,n}(\pi)-C_{k,n}(\pi)| & = & |\frac{1}{k^n} {n+k-d-1
\choose n} - \frac{1}{nk^{n-1}} {n+k-d-2 \choose n-1}|\\ & = &
|\frac{(n-d-1) (n+k-d-2)!}{k^n n! (k-d-1)!}|\\ & = & \frac{(n-d-1)
(n+k-d-2)!}{k^n n! (k-d-1)!}. \end{eqnarray*} Similarly if
$cd(\pi)=d(\pi)$, it follows that \begin{eqnarray*} | R_{k,n}(\pi)-C_{k,n}(\pi)|& = &
|\frac{1}{k^n} {n+k-d-1 \choose n} - \frac{1}{nk^{n-1}} {n+k-d-1
\choose n-1}|\\ & = & |\frac{-d (n+k-d-1)!}{k^n n! (k-d)!}|\\ & = &
\frac{d (n+k-d-1)!}{k^n n! (k-d)!}. \end{eqnarray*} 

	Letting $j$ index the number of cyclic descents, we conclude from Lemma \ref{cyclicdes} that \begin{eqnarray*}
& & ||R_{k,n}-C_{k,n}|| \\ & = & \frac{1}{2} \sum_{j=1}^{n-1} j A_{n-1,j} \frac{(n-j) (n+k-j-1)!}{k^n n! (k-j)!}\\ & & + (n-j) A_{n-1,j} \frac{j (n+k-j-1)!}{k^n n! (k-j)!} \\
& = & \frac{1}{k^n n} \sum_{j=1}^{n-1} j(n-j) A_{n-1,j} {n+k-j-1 \choose n-1}. \end{eqnarray*} This proves the first equality of the theorem.

	For the second equality of the theorem, consider the expression for $||R_{k,n}-C_{k,n}||$ given in the first equality of the theorem. We use the first part as follows. We break it into two terms as \[ \frac{1}{k^n} \sum_{j=1}^{n-1} j A_{n-1,j} {n+k-j-1 \choose n-1} - \frac{1}{nk^n} \sum_{j=1}^{n-1} j^2 A_{n-1,j} {n+k-j-1 \choose n-1}.\] To analyze the first term, we begin with part 3 of Lemma \ref{Worpitzky} (Worpitzky's identity): \[  \sum_{j=1}^{n-1} \frac{t^j A_{n-1,j}}{(1-t)^n} = \sum_{k=0}^{\infty} k^{n-1}t^k.\] Thus \[ \sum_{j=1}^{n-1} t^j A_{n-1,j} = (1-t)^n \sum_{k=0}^{\infty} k^{n-1}t^k.\] Differentiate both sides with respect to $t$. One obtains the equality \[ \sum_{j=1}^{n-1} jt^{j-1} A_{n-1,j} = -n(1-t)^{n-1} \sum_{k=0}^{\infty} k^{n-1}t^{k} + (1-t)^n \sum_{k=0}^{\infty} k^nt^{k-1}.\] Dividing both sides by $(1-t)^n$ gives \[ \sum_{j=1}^{n-1} \frac{jt^{j-1} A_{n-1,j}}{(1-t)^n} = \frac{-n}{1-t} \sum_{k=0}^{\infty} k^{n-1}t^k + \sum_{k=0}^{\infty} k^nt^{k-1}.\] Taking the coefficient of $t^{k-1}$ on both sides now gives the equation \[ \sum_{j=1}^{n-1} j A_{n-1,j} {n+k-j-1 \choose n-1} = k^n - n \sum_{s=0}^{k-1} s^{n-1}.\]

	Next consider the term \[ \frac{1}{nk^n} \sum_{j=1}^{n-1} j^2
A_{n-1,j} {n+k-j-1 \choose n-1}.\] Recall from the previous paragraph
that \[ \sum_{j=1}^{n-1} \frac{jt^{j-1} A_{n-1,j}}{n} = -(1-t)^{n-1}
\sum_{k=0}^{\infty} k^{n-1}t^{k} + \frac{(1-t)^n}{n}
\sum_{k=0}^{\infty} k^nt^{k-1}.\] Multiplying both sides by $t$ and
differentiating yields the equation \begin{eqnarray*} & &
\sum_{j=1}^{n-1} \frac{j^2 t^{j-1} A_{n-1,j}}{n}\\ & = & -(1-t)^{n-1}
\sum_{k \geq 0} (k+1) k^{n-1} t^k + (n-1)(1-t)^{n-2} \sum_{k \geq 0}
k^{n-1}t^{k+1}\\ & & -(1-t)^{n-1} \sum_{k \geq 0} k^nt^k +
\frac{(1-t)^n}{n} \sum_{k \geq 0} k^{n+1} t^{k-1}. \end{eqnarray*}
Dividing both sides by $(1-t)^n$ and then taking the coefficient of
$t^{k-1}$ one obtains that \begin{eqnarray*} & & \frac{1}{nk^n}
\sum_{j=1}^{n-1} j^2 A_{n-1,j} {n+k-j-1 \choose n-1}\\ & = &
\frac{k^{n+1}}{n} - (n+1) \sum_{s=1}^{k-2} s^n -
2(k-1)^n-\sum_{s=1}^{k-1} s^{n-1}\\
& & +(n-1)(k-1) \sum_{s=1}^{k-2}
s^{n-1}.\end{eqnarray*} This simplifies to the expression in the
second part of the theorem.

	The third expression follows from the second expression by the
	expansion of partial power sums in \cite{GR}: \[
	\sum_{r=0}^{a-1} r^n = \frac{a^n}{n+1} \left( a+
	\sum_{t=0}^{n-1} \frac{B_{t+1} (n+1)_{t+1}}{(t+1)! a^t} \right).\]
\end{proof}

	Theorem \ref{boundtot} will give upper and lower bounds on the total variation distance between $R_{k,n}$ and $C_{k,n}$, and implies that $log_2(n)+c$ iterations of a cut followed by a 2-riffle shuffle are necessary and sufficient for randomness. For the lower bound we will use the following elementary lemma about Bernoulli numbers.

\begin{lemma} \label{Bernoulli} For $t \geq 1$ an integer, the function $\frac{(-1)^{t-1} B_{2t}}{(2t-1)!}$ is a positive, decreasing function of $t$. \end{lemma}

\begin{proof} It is well known that $\frac{B_{2t}}{(2t)!} = \frac{(-1)^{t-1} 2 \zeta(2t)}{(2 \pi)^{2t}}$ where $\zeta$ is the Riemann zeta function. The result now follows from the fact that the Riemann zeta function is decreasing in $t$ together with the inequality $(2 \pi)^2 \geq \frac{t+1}{t}$ for $t \geq 1$. \end{proof}

\begin{theorem} \label{boundtot}
\begin{enumerate}
\item \[ ||R_{k,n}-C_{k,n}|| \leq .25 \frac{n}{k}. \] In particular, for $k \geq \frac{ne^{c_0}}{4}$, the total variation distance is at most $e^{-c_0}$.
\item \[ ||R_{k,n}-C_{k,n}|| \geq .24 \frac{n-1}{k} + \frac{(n+1)B_n}{k^{n-1}} \] for $k \geq n$. The term $\frac{(n+1)B_n}{k^{n-1}}$ goes to $0$ as $n \rightarrow \infty$ (again for $k \geq n$). 
\end{enumerate}
\end{theorem}

\begin{proof} For this upper bound, we use the first equality of Theorem \ref{exactformula}. Together with the fact that $j(n-j) \leq \frac{n^2}{4}$ for $1 \leq j \leq n-1$, it follows that \[ ||R_{k,n}-C_{k,n}|| \leq  \frac{n}{4k^n} \sum_{j=1}^{n-1} A_{n-1,j} {n+k-j-1 \choose n-1}.\] By part 3 of Lemma \ref{Worpitzky}, this is $\frac{n}{4k}$. 

	For the lower bound, we use the third equality of
	Theorem \ref{exactformula}, and suppose throughout that $k \geq n$. Then by Lemma \ref{Bernoulli} \[
	\sum_{t=1}^{n-2} \frac{B_{t+1}(n)_t}{t! k^t} \] is (after disregarding the terms corresponding to even $t$ which all vanish) an
	alternating sum whose terms decrease in magnitude, so is at
	least $\frac{n}{k} (B_2+\frac{B_4}{3!}) \geq .16
	\frac{n}{k}$. Similarly one sees that \[ \sum_{t=1}^{n-2}
	\frac{B_{t+1}(n-1)_t}{(t+1)!k^t} \geq \frac{n-1}{k} (\frac{B_2}{2}+\frac{B_4}{4!}) \geq .08 \frac{n-1}{k}.\] The assertion about the term involving $B_n$ follows from the well
	known asymptotics (\cite{O}) $B_{2t} \sim (-1)^t 4 \sqrt{\pi
	n} \left( \frac{n}{\pi e} \right)^{2n}$, so that for $k \geq n$ this term is extremely small. \end{proof}

	Next we record three corollaries of Theorem
	\ref{boundtot}. Corollary \ref{cutinc} relates to the
	distribution of the longest increasing subsequence.

\begin{cor} \label{cutinc} $log_2(n)+c$ iterations of a cut followed by a 2-riffle shuffle are sufficient to bring the distribution of the longest increasing subsequence close to that of a random permutation. \end{cor}

\begin{proof} This is immediate from Theorem \ref{boundtot} and the fact (see the discussion in Section \ref{shuffleinc} of this paper), that $\frac{5}{6}log_2(n)$ iterations of a 2-riffle shuffle are sufficient to bring the distribution of the longest increasing subsequence close to that of a random permutation. \end{proof}

	Corollary \ref{descentinversions} pertains to the numbers of
	descents.

\begin{cor} \label{descentinversions} $\frac{3}{2}log_2(n)+c$ iterations of a cut followed by a 2-riffle shuffle are necessary and sufficient to randomize the number of descents of a permutation.
\end{cor}

\begin{proof} It is known \cite{BD} that $\frac{3}{2}log_2(n)+c$ iterations of a 2-riffle shuffle are necessary and sufficient to random the number of descents of a permutation. The result now follows from Theorem \ref{boundtot}. \end{proof}

	Corollaries \ref{cutinc} and \ref{descentinversions} applied Theorem \ref{boundtot} to deduce information about the measure $C_{k,n}$ using information about the measure $R_{k,n}$. Corollary \ref{singlecard} gives an example where Theorem \ref{boundtot} is used to deduce information about the measure $R_{k,n}$ from information about the measure $C_{k,n}$. In fact the result of Corollary \ref{singlecard} is known \cite{U}, by an entirely different method of proof.

\begin{cor} \label{singlecard} $log_2(n)+c$ iterations of a 2-riffle shuffle are sufficient to randomize the position of any single card. \end{cor}

\begin{proof} Under the measure $C_{k,n}$, the position of any card is perfectly random. The result follows from Theorem \ref{boundtot}. \end{proof}  

\section{Affine Shuffles} \label{affine}

	This section considers a measure $A_{k,n}$ on the symmetric
        group known as the affine k-shuffle measure. It was introduced
        in \cite{F2} and studied further in \cite{F3}. As above, let
        $cd(\pi)$ be the number of cyclic descents of $\pi$. Let
        $maj(\pi)$ be the major index of $\pi$ (i.e. the sum of the values
        $i$ with $1 \leq i \leq n-1$ such that $\pi(i)>\pi(i+1)$). Let
        $C_r(m)$ denote the Ramanujan sum $\sum_{1 \leq l \leq r \atop
        gcd(l,r)=1} e^{2 \pi i l m/r}$. Then the measure $A_{k,n}$
        assigns mass \[ \begin{array}{ll}
        \frac{1}{nk^{n-1}} \sum_{r|n,k-cd(\pi)}\
        {\frac{n+k-cd(\pi)-r}{r} \choose \frac{n-r }{r}}
        C_r(-maj(\pi))& \mbox {if \ $k-cd(\pi) >0$}\\ \
        \frac{1}{k^{n-1}} & \mbox{if \ $k=cd(\pi), n | maj(\pi)$} \\
        0 & \mbox{otherwise}. \end{array} \] Note that the $r=1$ term in this formula is the measure $C_{k,n}$, suggesting that $A_{k,n}$ and $C_{k,n}$ should be quite close.

	When $k=2$ and the deck size $2n$ is even, the inverse of an affine shuffle
has a physical description: choose an even number $2j$ between
$1$ and $2n$ with the probability of getting $2j$ equal to $\frac{{2n
\choose 2j}}{2^{2n-1}}$. From the stack of $2n$ cards, form a second
pile of size $2j$ by removing the top $j$ cards of the stack and
putting the bottom $j$ cards of the first stack on top of them; then riffle the two stacks (of size $2n-2j$ and $2j$) together as in ordinary riffle shuffles. One of the main results of this section is that as for ordinary riffle shuffles, $\frac{3}{2}log_2(n)+c$ iterations of an affine 2-shuffle are necessary and sufficient for randomness.

	The probability measure $A_{k,n}$ has many remarkable properties. First, it satisfies the equation \[ \left( \sum_{\pi \in S_n} A_{k_1,n}(\pi) \pi \right) \left( \sum_{\pi \in S_n} A_{k_2,n}(\pi) \pi \right) = \left( \sum_{\pi \in S_n} A_{k_1k_2,n}(\pi) \pi \right). \] Second, one can understand in a very precise way the distribution of this measure on conjugacy classes, even though the measure $A_{k,n}$ is not constant on conjugacy classes. See \cite{F2},\cite{F3} for details. 

	Theorem \ref{affineclosecut} gives an upper bound on the total variation distance $||A_{k,n}-C_{k,n}||$, which is sufficient for our purposes. 

\begin{theorem} \label{affineclosecut} $||A_{k,n}-C_{k,n}||< n \left(\frac{2n}{ek} \right)^{n/2}$ for $k \geq n$. \end{theorem} 

\begin{proof} Since $k \geq n$ and $cd(\pi) \leq n-1$, clearly $k-cd(\pi)>0$ for all $\pi$. From the explicit formulas for $A_{k,n}$ and $C_{k,n}$ mentioned in this paper, it follows that \[||A_{k,n}-C_{k,n}|| = \frac{1}{2nk^{n-1}} \sum_{\pi} | \sum_{r|n,k-cd(\pi) \atop r>1} {\frac{n+k-cd(\pi)}{r}-1 \choose \frac{n}{r}-1} C_r(-maj(\pi))|.\] Using the fact that there are $n!$ permutation on $n$ symbols and the bound $|C_r| \leq \phi(r)$ where $\phi$ is Euler's $\phi$ function of elementary number theory, this is at most \[ \frac{n!}{2nk^{n-1}} \sum_{r|n \atop r>1} {\frac{n+k}{r} \choose \frac{n}{r}-1} \phi(r).\] Since $k \geq n$, \[ {\frac{n+k}{r} \choose \frac{n}{r}-1} \leq  \frac{( \frac{n+k}{r})^{\frac{n}{r}-1}}{(\frac{n}{r}-1)!}
\leq \frac{( \frac{2k}{r})^{\frac{n}{r}-1}}{(\frac{n}{r}-1)!} = \frac{n}{2k} \frac{( \frac{2k}{r})^{\frac{n}{r}}}{(\frac{n}{r})!}.\] Since $r \geq 2$ and $k \geq n$, this is less than $\frac{n k^{n/2}}{2k (n/2)!}$. Thus \[ ||A_{k,n}-C_{k,n}|| \leq \frac{n!}{4k^{n/2} (\frac{n}{2})!} \sum_{r|n,r>1} \phi(r).\] Using the fact that $\sum_{r|n} \phi(r)=n$, together with the Stirling formula bounds proved in \cite{Fe}
\[ n! < (2 \pi)^{1/2}n^{n+1/2} e^{-n+1/(12n)}\]
\[ (n/2)! > (2 \pi)^{1/2}(n/2)^{n/2+1/2} e^{-n/2+1/(6n+1)},\] it follows that \[ ||A_{k,n}-C_{k,n}|| \leq n \left( \frac{2n}{ek} \right)^{n/2}. \] \end{proof}

	Corollary \ref{easycor} is immediate from Theorem \ref{affineclosecut} and the corresponding results for the measure $C_{k,n}$ in Section \ref{totvarcut}.

\begin{cor} \label{easycor}
\begin{enumerate}
\item $\frac{3}{2}log_2(n)+c$ affine 2 shuffles are necessary and suffice for randomness.

\item $log_2(n)+c$ iterations of an affine 2 shuffle suffice to bring the longest increasing subsequence close to that of a random permutation.

\item $\frac{3}{2}log_2(n)+c$ iterations of an affine 2 shuffle are necessary and sufficient to randomize the number of descents of a permutation.

\item $log_2(n)+c$ iterations of an affine 2 shuffle are sufficient to randomize the position of any single card.
\end{enumerate}
\end{cor}

\section{Acknowledgements} The author was support by NSA Grant MDA904-03-1-0049. We thank Ira Gessel for discussions.

\end{document}